\input amstex
\input amsppt.sty

\input label.def
\input degt.def

\copycounter\thm\subsection
\newtheorem{Remark}\Remark\thm\endAmSdef


{\catcode`\@=11
\gdef\proclaimfont@{\sl}}

\loadbold

\def\bA{\bold A}
\def\bD{\bold D}
\def\bE{\bold E}
\def\bJ{\bold J}

\def\bX{\bold X}

\def\FF{\Bbb F}
\def\DD{\Bbb D}
\def\PP{\Bbb P}

\def\GB{\Bbb B} 
\def\GS{\Bbb S} 
\def\GD{\DD} 

\let\CK\CalK

\let\CO\CalO

\def\tB{\tilde B}
\def\tS{\tilde S}
\def\tX{\tilde X}

\def\KK{K}
\def\BK{\bar K}

\let\Ga\alpha
\let\Gb\beta
\let\Gf\varphi

\def\barO{\bar O}
\def\barD{\bar D}
\def\barP{\bar P}
\def\barL{\bar L}

\def\Tor{\operatorname{Tor}}

\def\Cp#1{\PP^{#1}}
\def\<#1>{\langle#1\rangle}
\def\ls|#1|{\mathopen|#1\mathclose|}

\def\ie,{\emph{i.e.},}
\def\eg,{\emph{e.g.},}
\def\etc.{\emph{etc}.}
\def\cf.{\emph{cf}\.}

\let\splus\oplus 
\let\HS=Z        

\def\onto{\to\joinrel\joinrel\joinrel\joinrel\to}

\topmatter

\author
Alex Degtyarev
\endauthor

\title
Oka's conjecture on irreducible plane sextics. II
\endtitle

\address
Department of Mathematics,
Bilkent University,
06800 Ankara, Turkey
\endaddress

\email
degt\@fen.bilkent.edu.tr
\endemail

\abstract
We complete the proof of Oka's conjecture on the
Alexander polynomial
of an irreducible plane sextic.
We also calculate the fundamental
groups of irreducible sextics with a singular point
adjacent to~$\bJ_{10}$.
\endabstract

\keywords
Plane sextic, torus type, Alexander polynomial,
dihedral covering, fundamental group
\endkeywords

\subjclassyear{2000}
\subjclass
Primary: 14H30; 
Secondary: 14J27 
\endsubjclass

\endtopmatter

\document

\section{Introduction}\label{S.intro}

This paper is a sequel to my paper~\cite{degt.Oka}, where Oka's
conjecture on the Alexander polynomial is settled for all
irreducible sextics with simple singularities. Here, we complete
the proof for the missing case of sextics with a non-simple
singular point adjacent to~$\bJ_{10}$ (a point of simple tangency
of three smooth branches).

Recall that a plane sextic $C\subset\Cp2$ is said to be of
\emph{torus type} if its equation can be represented in the form
$p^3+q^2=0$, where $p$ and~$q$ are some homogeneous polynomials of
degree~$2$ and~$3$, respectively (see Section~\ref{s.torus} for
details). Sextics of torus type are a
major source of examples of plane curves with large fundamental
group $\pi_1(\Cp2\sminus C)$. Historically, it was a six cuspidal
sextic of torus type that was the first example of an irreducible
plane curve with infinite fundamental group, see
O.~Zariski~\cite{Zariski}. (An irreducible quintic with infinite
fundamental group was discovered much later in~\cite{groups}.)
All sextics of torus type have
nontrivial Alexander polynomials (see below); hence, their
fundamental groups are infinite.

The fundamental group $\pi_1(\Cp2\sminus C)$ is a powerful
invariant of a plane curve, but it is extremely difficult to
calculate. A much simpler invariant, capturing the abelinization
of the first commutant of the group, is the so called
\emph{Alexander polynomial}. For an irreducible
curve $C\subset\Cp2$ of degree~$m$, its Alexander polynomial
$\Delta_C(t)$ can be defined as the characteristic
polynomial of the deck translation automorphism of the vector
space $H_1(X_m\sminus C;\C)$, where $X_m\to\Cp2$ is the
cyclic $m$-fold covering ramified at~$C$. (Since $C$ is assumed
irreducible, such a covering is unique.)
The Alexander polynomial
is a purely algebraic invariant of the fundamental group;
in particular, $\Delta_C(t)\ne1$
if and only if the
quotient $K/[K,K]$ is infinite, where $K$ is the commutant of
the group.
For more details, alternative definitions, and basic properties of the
Alexander polynomial see A.~Libgober~\cite{Libgober}
and~\cite{Libgober2}.
The Alexander polynomials of irreducible sextics are calculated
in~\cite{poly}. For further references, see M.~Oka's
survey~\cite{Oka.survey}.

Note that the Alexander polynomial of a plane curve is subject to
rather strong divisibility conditions, see~\cite{Zariski},
\cite{Libgober}, and~\cite{divide}. In particular, six is the
first degree where the polynomial of an irreducible curve may be
nontrivial.

Based on the known examples, Oka~\cite{Oka.conjecture} conjectured
that {\proclaimfont any irreducible sextic whose Alexander
polynomial is nontrivial is of torus type\/}. (A similar
conjecture on the fundamental group was disproved
in~\cite{degt.Oka}; more counterexamples are given by
Theorem~\ref{special} below.) Proof of Oka's conjecture
constitutes the main result of the present paper. More precisely,
the following statement holds.

\theorem[Main Theorem]\label{th.main}
For an irreducible plane sextic~$C$, the
following statements are equivalent\rom:
\roster
\item\local{torus}
$C$ is of torus type\rom;
\item\local{Delta}
the Alexander polynomial $\Delta_C(t)$ is nontrivial\rom;
\item\local{B3}
the group $\pi_1(\Cp2\sminus C)$ factors to the
reduced braid group~$\GB_3/\Delta^2$\rom;
\item\local{S3}
the group $\pi_1(\Cp2\sminus C)$ factors to the
symmetric group~$\GS_3$.
\endroster
\endtheorem

Theorem~\ref{th.main} is proved
in~\cite{degt.Oka} under the assumption that either all singular
points of~$C$ are simple or $C$ has a non-simple singular point
adjacent to~$\bX_9$ (a quadruple point). It is also shown
in~\cite{degt.Oka} that \loccit{torus} implies
\loccit{Delta}--\loccit{S3} and that \loccit{Delta}
implies~\loccit{S3}. (Obviously, \loccit{B3} implies~\loccit{S3}
as well.)
The only remaining case is the implication
\loccit{S3}~$\Rightarrow$~\loccit{torus} for
a sextic with a
triple non-simple singular point. All such points are adjacent
to~$\bJ_{10}$ (a semiquasihomogeneous singularity of type $(6,3)$); a
sextic with such a point is called a \emph{$\bJ$-sextic}, see
Definition~\ref{def.J}.
Thus, Theorem~\ref{th.main} is a consequence of
the following statement, which is actually proved in the paper.

\theorem\label{main.J10}
Let~$C$ be an irreducible $\bJ$-sextic,
and assume that the fundamental group
$\pi_1(\Cp2\sminus C)$ factors to~$\GS_3=\GD_6$.
Then $C$ is of torus type.
\endtheorem

In fact, Theorem~\ref{main.J10} holds for reducible $\bJ$-sextics
as well, see Remark~\ref{rem.torus}.

As in the case of sextics with simple singularities,
see~\cite{degt.Oka},
Theorem~\ref{main.J10} admits a slightly more precise version
(which is meaningful for the sets of singularities
$\bJ_{2,0}\oplus4\bA_2$ and $\bJ_{2,3}\oplus3\bA_2$).

\theorem\label{main.J10'}
For an irreducible $\bJ$-sextic~$C$, there is a natural one to
one correspondence between the set of quotients
of $\pi_1(\Cp2\sminus C)$ isomorphic to~$\GD_6$
and the set of torus structures of~$C$.
\endtheorem

Throughout the paper, we use the notation of~\cite{Arnold} for the
types of singularities adjacent to $\bJ_{10}=\bJ_{2,0}$. One can
use Table~\ref{table1} below as a guide. Although, in the presence
of non-simple singular points, a set of singularities is no longer
determined by its resolution lattice, we keep using the lattice
notation~$\splus$ in the listings. (In general, we refer
to~\cite{degt.Oka} for the less common notation and terminology.)

As another application of the approach developed in the paper, we
calculate the fundamental groups of most irreducible
$\bJ$-sextics (all except the two families mentioned
in Theorem~\ref{torus}).
The results are stated in
Theorems~\ref{special}--\ref{abelian} below.

\theorem\label{special}
There exist irreducible plane sextics with the
following
sets of singularities\rom: $\bJ_{2,0}\splus2\bA_4$,
$\bJ_{2,1}\splus2\bA_4$, and $\bJ_{2,5}\splus\bA_4$.
For each such sextic~$C$,
one has $\pi_1(\Cp2\sminus C)=\GD_{10}\times(\ZZ/3\ZZ)$, where
$\GD_{10}$ is the dihedral group of order~$10$.
\endtheorem

\theorem\label{torus}
Let~$C$ be an irreducible $\bJ$-sextic of torus type
and with a set of singularities
other than $\bJ_{2,0}\splus4\bA_2$ or $\bJ_{2,3}\splus3\bA_2$.
Then
$\pi_1(\Cp2\sminus C)$ is
the reduced braid group $\GB_3/\Delta^2$.
\endtheorem

For the two exceptional sets of singularities
listed in Theorem~\ref{torus}, the
Alexander polynomial of the curve is $(t^2-t+1)^2$, see,
\eg,~\cite{poly}. Hence, the fundamental group must be much larger
than~$\GB_3/\Delta^2$. As in Section~\ref{proof.torus}, one
can
use~\cite{groups} and
assert that this group is a quotient of
$$
\bigl<a,b,c\bigm|
 aba=bab,\ bcb=cbc,\ abcb^{-1}a=bcb^{-1}abcb^{-1}\bigr>.
$$

\theorem\label{abelian}
Let~$C$ be an irreducible $\bJ$-sextic that is neither of torus
type nor one of the curves listed in Theorem~\ref{special}.
Then the group $\pi_1(\Cp2\sminus C)$ is
abelian.
\endtheorem

\subsection{Contents of the paper}
Proof of Theorem~\ref{main.J10} follows the lines of~\cite{poly}:
we construct an appropriate conic $Q=\{p=0\}$ and cubic
$K=\{q=0\}$ and, using the B\'ezout theorem, show that the
difference $\Gf-q^2$ should be divisible by~$p^3$ (where
$\{\Gf=0\}$ is the original curve).

To simplify the calculation
and minimize the number of cases to be considered, we start with
reducing the problem to the study of the so called \emph{trigonal
models}, which are trigonal curves on the quadratic cone
in~$\Cp3$. This is done in~\S\ref{S.J10}. We show that, in all
cases of interest, trigonal models have simple singularities (and
in rather small number: the set of singularities should admit an
embedding to~$\bE_8$)
and extend to such curves the results of~\cite{JAG}
and~\cite{degt.Oka} concerning plane sextic. An important result
here is Proposition~\ref{B=C}, which relates the fundamental
groups of a $\bJ$-sextic and its trigonal model, and its
Corollary~\ref{cor.B=C} stating that, to a certain extent, the
fundamental group of a $\bJ$-sextic does not depend on its
non-simple singular point.

Theorems~\ref{main.J10} and~\ref{main.J10'} are
proved in \S\ref{S.proof}. Again, the
torus type condition for the original curve is reduced to a
similar condition for its trigonal model, so that one needs to
deal with simple singularities only. Note that, unlike the case of
abundant \emph{vs\.} non-abundant curves in~\cite{poly}, where the
conic~$Q$ was obtained from the linear system calculating the
Alexander polynomial of the curve, here the existence of~$Q$
follows from a simple dimension count.

In~\S\ref{S.groups}, we continue exploring properties of
trigonal models and find the fundamental groups of most
$\bJ$-sextics. Theorems~\ref{special}--\ref{abelian} are proved
here. We use a simple consequence of van Kampen's
method~\cite{vanKampen} (avoiding any attempt to calculate the
global braid monodromy) and reduce the few remaining cases to the
results of~\cite{groups} dealing with plane quintics.

In concluding \S\ref{S.application}, we apply the results of
previous sections to prohibit several $\bJ$-sextics.
Propositions~\ref{no.J10} and~\ref{no.E12} proved here are not
new, see~\cite{quartics}, but we give new short proofs. It is
remarkable that prohibited are precisely the curves that would
contradict Theorems~\ref{main.J10}--\ref{abelian}.

\section{The trigonal model}\label{S.J10}

In this section, we introduce the so called \emph{trigonal model}
of a $\bJ$-sextic~$C$, thus reducing the study of~$C$ to the study
of a certain curve with simple singularities.

\definition\label{def.J}
A \emph{$\bJ$-sextic} is a reduced plane sextic $C\in\Cp2$ with a
distinguished
triple singular point~$O$ adjacent to~$\bJ_{10}$ and without linear
components passing through~$O$.
\enddefinition

Let $(C,O)$ be a $\bJ$-sextic. Blow up~$O$ to obtain a Hirzebruch
surface~$\Sigma_1=\Cp2(O)$. The proper transform of~$C$ has a
triple point~$O'$ in the exceptional section. Blow it up and blow
down the fiber through~$O'$. The result is a Hirzebruch
surface~$\Sigma_2$ (quadratic cone with the vertex blown up); we
denote it by~$\HS$. It is a geometrically ruled surface with an
exceptional section~$S_0$ of selfintersection~$(-2)$. When
speaking about \emph{fibers} of~$\HS$, we mean fibers of the
ruling.

The proper transform of~$C$ is a certain curve
$B\subset\HS$ disjoint from~$S_0$;
it is called the \emph{trigonal model} of~$C$. The
trigonal model is equipped with a distinguished point $\barO$
(image of the fiber through~$O'$) and distinguished fiber
$F_0\ni\barO$ (image of the exceptional divisor over~$O'$). It is
important that $\barO$ does not belong to either~$B$ or~$S_0$.
Below, we show that, as long as the fundamental group is
concerned, $\barO$ and~$F_0$ are irrelevant (see
Proposition~\ref{B=C}).

Let~$F$ be a generic fiber of~$\HS$. As is well known, $S_0$
and~$F$ generate the semigroup of effective divisors on~$\HS$. One
has $B\in\ls|3S_0+6F|$; in particular, $B$ \emph{is} a trigonal
curve, in the sense that it intersects each fiber of~$\HS$ at
three points. The \emph{singular fibers} of~$B$ are the fibers
of~$\HS$ that are not transversal to~$B$. Counting with
multiplicities, $B$ has twelve singular fibers.
If~$B$ has simple singularities only, its singular fibers can be
regarded as the singular fibers of the Jacobian elliptic
surface~$\tX$ obtained as the double covering of~$\HS$ ramified
at~$B$ and~$S_0$, see below.

Note that $B$ may be reducible; however, it cannot contain~$S_0$
or a fiber of~$\HS$. Any curve $B\in\ls|3S_0+6F|$ satisfying this
condition
is called a
\emph{trigonal model}. Any trigonal model~$B$ gives rise to a
$\bJ$-sextic after a point $\barO\in\HS\sminus(B\cup S_0)$ is
chosen.

\lemma\label{simple}
If a $\bJ$-sextic~$C$ is irreducible, then all singular points of
its trigonal model~$B$ are simple.
\endlemma

\proof
The statement follows immediately from the genus formula, as a
nonsingular curve in $\ls|3S_0+6F|$ has genus~$4$ and any
non-simple singular point takes off the genus at least~$6$.
\endproof

It is easy to see (\eg, using the associated cubic, see
Section~\ref{s.cubic}
below) that the only non-simple singular point that $B$ may have
is~$\bJ_{10}$. In this case, depending on the distinguished
fiber~$F_0$, the set of singularities of the original sextic~$C$
is either $2\bJ_{10}$ or~$\bJ_{4,0}$, and $C$ consists of three
conics, which either have two common tangency points (the
case~$2\bJ_{10}$) or one common point of $4$-fold intersection
(the case~$\bJ_{4,0}$). Note that both families are obviously of
torus type.

From now on, we always assume that all singular points of~$B$ are
simple. We identify a set of simple singularities of~$B$ with its
resolution lattice, which is a direct sum of irreducible root
systems ($\bA$--$\bD$--$\bE$ lattices), one summand for each
singular point of the same name.
The relation between the types of the distinguished
fiber~$F_0$ and distinguished singular point~$O$ of~$C$ is given by
Table~\ref{table1}, where for~$F_0$ we use (one of) the standard
notation for singular elliptic fibers.

\midinsert
\table
Singular fibers and singular points
\endtable\label{table1}
\vskip-\bigskipamount
\line{\offinterlineskip
\def\header{\vtop\bgroup
 \halign\bgroup\vrule height12pt depth3.5pt width0pt
  \vrule\quad$##$\quad\hss\vrule&\quad$##$\quad\hss\vrule\cr
 \noalign{\hrule}
 \omit\vrule\vrule height11pt depth3.5pt width0pt\hss\ Fiber~$F_0$\ \hss\vrule
  &\omit\hss\ Point~$O$\ \hss\vrule\cr
 \noalign{\hrule}}
\def\footer{\noalign{\hrule}\crcr\egroup\egroup}
\hss\header
\tilde\bA_0&\bJ_{2,0}\cr
\tilde\bA_0^*&\bJ_{2,1}\cr
\tilde\bA_0^{**}&\bE_{12}\cr
\tilde\bA_1^*&\bE_{13}\cr
\footer\quad
\header
\tilde\bA_2^*&\bE_{14}\cr
\tilde\bE_6&\bE_{18}\cr
\tilde\bE_7&\bE_{19}\cr
\tilde\bE_8&\bE_{20}\cr
\footer\quad
\header
\tilde\bA_p$, $p\ge1&\bJ_{2,p+1}\cr
\tilde\bD_q$, $q\ge4&\bJ_{3,q-4}\cr
\footer\hss}\null
\endinsert

\proposition\label{B=C}
For a $\bJ$-sextic~$(C,O)$ and its trigonal model $B\subset\HS$,
there is a canonical isomorphism
$\pi_1(\Cp2\sminus C)=\pi_1(\HS\sminus(B\cup S_0))$.
\endproposition

\proof
Consider the surface $\tilde\HS$ obtained from~$\HS$ by blowing up
the distinguished point~$\barO$. Denote by~$\sptilde$ the proper
pull-backs of the curves involved, and let $\tilde L$ be the
exceptional divisor over~$\barO$. There are obvious
diffeomorphisms
$$
\gather
\Cp2\sminus C
 =\tilde\HS\sminus(\tilde B\cup\tilde S_0\cup\tilde F_0),\\
\Cp2\sminus(C\cup L)
 =\tilde\HS\sminus(\tilde B\cup\tilde S_0\cup\tilde F_0\cup\tilde L)
 =\HS\sminus(B\cup S_0\cup F_0),
\endgather
$$
where $L\subset\Cp2$ is the line tangent to~$C$ at~$O$. Thus,
$\pi_1(\Cp2\sminus C)$ is obtained from the group
$\pi_1(\tilde\HS\sminus(\tilde B\cup\tilde S_0\cup\tilde F_0\cup\tilde L))$
by adding the relation $[\partial\tilde\Gamma]=1$, where
$\tilde\Gamma\subset\tilde\HS$ is a small analytic disc
transversal to~$\tilde L$ (and disjoint from the other curves).
The projection $\Gamma\subset\HS$ of~$\tilde\Gamma$ is a small
analytic disk transversal to~$F_0$ at~$\barO$; since $\barO$ does
not belong to the union $B\cup S_0$, the relation
$[\partial\Gamma]=1$ is precisely the relation resulting from
patching the distinguished fiber~$F_0$.
\endproof

\corollary\label{cor.B=C}
The fundamental group of a $\bJ$-sextic~$C$ obtained from a
trigonal model~$B$ does not depend on the choice of a
distinguished fiber~$F_0$.
\qed
\endcorollary

Let $p\:X\to\HS$ be the double covering of~$\HS$ ramified along
$B+S_0$, and let $\tX$ be the minimal resolution of singularities
of~$X$. It is well known that $\tX$ is a rational elliptic
surface; its intersection lattice $H_2(\tX)$ is the only odd
unimodular lattice of signature $(1,9)$. Let~$E$ be the
exceptional divisor contracted by the projection $\tX\to X$, and
let~$\tS_0$ and~$\tB$ be
the proper transforms of, respectively, $S_0$ and~$B$ in~$\tX$.
The copies of~$S_0$ and~$B$ in~$X$ are identified with~$S_0$
and~$B$ themselves.

Let $s_0=[\tS_0]\in H_2(\tX)$, let $f\in H_2(\tX)$ be the
class realized by the pull-\hskip0ptback of a generic fiber
of~$\HS$, and denote by $T_B$ the sublattice spanned by~$s_0$
and~$f$, \ie, $T_B=\ZZ s_0+\ZZ f\subset H_2(\tX)$.

\lemma\label{E8}
The sublattice $T_B\subset H_2(\tX)$ is an orthogonal
summand, and one has
$T_B^\perp\cong\bE_8$.
\endlemma

\proof
One has $e_0^2=-1$, $f^2=0$, and $e_0\cdot f=1$. Hence, $T_B$ is a
unimodular lattice of signature~$(1,1)$. In particular, $T_B$ is
an orthogonal summand in any larger lattice. The orthogonal
complement~$T_B^\perp$ is a unimodular lattice of signature
$(0,8)$ and, to complete the proof, it remains to show that
$T_B^\perp$ is even, \ie, that $T_B$ contains a characteristic
vector of $H_2(\tX)$.

Replace~$\tX$ with a surface~$X'$ obtained from~$X$ by a small
perturbation; it can be regarded as the double covering of~$\HS$
ramified at~$S_0$ and a nonsingular curve~$B'$ obtained by a
perturbation of~$B$. There is a diffeomorphism $\tX\cong X'$
identical outside some regular neighborhoods of the singular points,
see~\cite{Durfee}. Since $\HS$ is a Spin-manifold, the topological
projection formula implies that the Stiefel-Whitney class
$w_2(X')$ is given by $([B']+[S_0])\bmod2\in H_2(X';\ZZ/2\ZZ)$. On
the other hand, using the fact that $H_2(X')$ is torsion free, one
can see that $[B']=3s_0+3f$. Hence, $w_2(\tX)=f\bmod2$, and the
statement follows.
\endproof

Let~$\Sigma_B$ be the set of singularities of~$B$. Recall that we
identify a set of simple singularities with its resolution
lattice, \ie, $\Sigma_B$ can be
regarded as the sublattice in $H_2(\tX)$ spanned by the classes of
the exceptional divisors in~$\tX$ (those that are contracted
in~$X$). Obviously, $\Sigma_B\subset T_B^\perp\cong\bE_8$. Denote
$\CK_B=\Tors(T_B^\perp/\Sigma_B)$.

Consider the canonical epimorphism
$\kappa\:\pi_1(\HS\sminus(B\cup S_0))\onto\ZZ/2\ZZ$. (If $B$ is
irreducible, then $\kappa$ is the only epimorphism to $\ZZ/2\ZZ$.
Otherwise, $\kappa$ can be defined as the map sending each
van Kampen generator
to $1\in\ZZ/2\ZZ$. Alternatively,
$\Ker\kappa\subset\pi_1(\HS\sminus(B\cup S_0))$ is the group of
the covering $X\sminus(B\cup S_0)\to\HS\sminus(B\cup S_0)$.)
Let $\KK(B)=\Ker\kappa$
and denote by $\BK(B)$ the abelinization of~$\KK(B)$.
One has
$$
\KK(B)=\pi_1(X\sminus(B\cup S_0))
\quad\text{and}\quad
\BK(B)=H_1(X\sminus(B\cup S_0)).
$$
Denote by~$\tr$ the automorphism of $\BK(B)$
given by $[a]\mapsto[\tilde1^{-1}a\tilde1]$, where $a\in\KK(B)$,
$[a]$ stands for the class realized by~$a$ in $\BK(B)$, and
$\tilde1\in\pi_1(\HS\sminus(B\cup S_0))$
is a lift of the generator $1\in\ZZ/2\ZZ$.
Alternatively, $\tr$
is induced by the deck translation of the covering
$X\sminus(B\cup S_0)\to\HS\sminus(B\cup S_0)$.

Next three statements are analogs of similar statements
for plane sextics with simple singularities,
see~\cite{JAG} and~\cite{degt.Oka}.
Proofs are omitted; instead, we refer to the counterparts for
plane sextics, whose proofs apply almost literally.

\proposition[Proposition \rm(\cf. Theorem~4.3.1 in~\cite{JAG})]\label{irreducible}
A trigonal model~$B$ \rom(and, hence, a $\bJ$-sextic~$C$\rom) is
irreducible if and only if the group~$\CK_B$ is free of
$2$-torsion.
\qed
\endproposition

\proposition[Proposition \rm(\cf. Proposition~3.4.4 in~\cite{degt.Oka})]\label{splitting}
Let~$B$ be
an irreducible trigonal model.
Then there is a splitting
$\BK(B)=\Ker(\tr-1)\oplus\Ker(\tr+1)$
and canonical isomorphisms
$\Ker(\tr-1)=\ZZ/3\ZZ$ and $\Ker(\tr+1)=\Ext(\CK_B,\ZZ)$.
\qed
\endproposition

\corollary[Corollary \rm(\cf. Corollary~3.4.6 in~\cite{degt.Oka})]\label{dihedral}
Let~$B$ be the trigonal model of an irreducible $\bJ$-sextic~$C$.
Then
there is a canonical one to one correspondence between the set of
normal subgroups $N\subset\pi_1(\Cp2\sminus C)$ with
$\pi_1(\Cp2\sminus C)/N\cong\GD_{2n}$, $n\ge3$,
and the set of subgroups of $\Tor(\CK_B,\ZZ/n\ZZ)$ isomorphic to
$\ZZ/n\ZZ$.
\qed
\endcorollary

Lemma~\ref{E8} and Proposition~\ref{irreducible} result in the
following necessary condition for a set of simple singularities to
be realized by the trigonal model of a $\bJ$-sextic.

\corollary\label{embedding}
Let~$\Sigma$ be a set of simple singularities of a trigonal
model~$B$. Then $\Sigma$ admits an embedding to~$\bE_8$.
If $B$ is irreducible, then there is an embedding
$\Sigma\hookrightarrow\bE_8$
with $\bE_8/\Sigma$ free of $2$-torsion.
\qed
\endcorollary

\section{Proof of Theorems~\ref{main.J10} and~\ref{main.J10'}}\label{S.proof}

As in the previous section, we fix a $\bJ$-sextic~$C$
and denote by $B\subset\HS$ its trigonal model. We assume that all
singular points of~$B$ are simple, see Lemma~\ref{simple}.


\lemma\label{Sigma.E8}
Let~$\Sigma$ be a root system embedded to~$\bE_8$,
and let $\CK=\Tors(\bE_8/\Sigma)$.
If $\ls|\CK|$ is odd \rom(in particular, $\CK\ne0$\rom),
then $\Sigma$ is one of the following lattices\rom:
\widestnumber\item{$\CK=(\ZZ/3\ZZ)^2:$}
\setbox0\hbox{$\CK=(\ZZ/3\ZZ)^2:$}
\dimen0\wd0
\def\bx#1{\hbox to\dimen0{$#1$\hss}}
\roster
\item"\bx{\CK=(\ZZ/3\ZZ)^2:}"
$\Sigma=4\bA_2$;
\item"\bx{\CK=\ZZ/3\ZZ:}"
$\Sigma=3\bA_2$, $3\bA_2\oplus\bA_1$, $\bA_5\oplus\bA_2$,
$\bA_8$, or $\bE_6\oplus\bA_2$;
\item"\bx{\CK=\ZZ/5\ZZ:}"
$\Sigma=2\bA_4$.
\endroster
Conversely, for each root system listed above there is a unique,
up to isomorphism, embedding to~$\bE_8$.
\endlemma

\proof
The enumeration of the embeddings of root systems to~$\bE_8$ is a
simple task. For example, one can use V.~Nikulin's
techniques~\cite{Nikulin} of lattice extensions and discriminant
forms. We omit the details.
\endproof

\corollary\label{D6}
The fundamental group of an irreducible $\bJ$-sextic~$C$ factors
to the dihedral group~$\GD_6$ \rom(respectively,~$\GD_{10}$\rom)
if and only if the set of
singularities of the trigonal model of~$C$ is one of the
following: $3\bA_2\splus\ldots$, $\bA_5\splus\bA_2$,
$\bA_8$, or $\bE_6\splus\bA_2$ \rom(respectively, $2\bA_4$\rom).
\endcorollary

\proof
The statement follows from
Corollary~\ref{dihedral} and Lemma~\ref{Sigma.E8}.
\endproof

\subsection{Torus structures}\label{s.torus}
A plane sextic~$C$ is said to be of \emph{torus type} if its
equation can be represented in the form $p^3+q^2=0$, where $p$
and~$q$ are some homogeneous polynomials in $(x_0,x_1,x_2)$
of degree~$2$ and~$3$,
respectively. Any representation as above (considered up to
rescaling) is called a \emph{torus structure} of~$C$. With the
exception of a few very degenerate cases, a torus structure is
determined by the conic $Q=\{p=0\}$.

A sextic is of torus type if and only if it is the critical locus
of a projection to~$\Cp2$ of a cubic surface $V\subset\Cp3$; the
latter is given by $3x_3^3+3x_3p+2q=0$.

Each point of intersection of the conic $Q=\{p=0\}$ and cubic
$K=\{q=0\}$ is a singular point of~$C$; such points are called
\emph{inner}, and the other singular points that $C$ may have are
called \emph{outer}. The type of a simple inner singular point~$P$
is determined by the mutual topology of~$Q$ and~$K$ at~$P$,
whereas outer points occur in the family $(\Ga p)^3+(\Gb q)^2=0$
under some special values of parameters $\Ga,\Gb\in\C^*$.
Note that, in the case of non-simple singularities, one can speak
about `outer degenerations' of inner singular points. Thus,
with~$Q$ and~$K$ fixed, an inner point of
type~$\bJ_{10}=\bJ_{2,0}$ may degenerate to $\bJ_{2,1}$ or
$\bJ_{2,2}$.

\lemma\label{reduction}
A $\bJ$-sextic~$C$ is of torus type if and only if there exist
sections
$p\in\Gamma(\HS;\CO_\HS(S_0+2F))$ and
$q\in\Gamma(\HS;\CO_\HS(S_0+3F))$ such that
the trigonal model~$B$ of~$C$ is given by an equation of the form
$p^3+s_0q^2=0$, where
$s_0\in\Gamma(\HS;\CO_\HS(S_0))$ is a fixed section whose zero set
is~$S_0$.
\endlemma

\proof
Let $\bar\Gf=0$ be an equation of~$C$, and let $\bar\Gf=\bar p^3+\bar q^2$
be its torus structure. It is well known that the
conic $Q=\{\bar p=0\}$ at~$O$ is smooth and tangent to~$C$,
the cubic $K=\{\bar q=0\}$ is singular at~$O$, and the local
intersection index of~$Q$ and~$K$ at~$O$ is at least~$3$. (Indeed,
if both~$K$ and~$Q$ are singular, then $C$ has a quadruple point
at~$O$. If $K$ is nonsingular or the local intersection index
of~$Q$ and~$K$ at~$O$ is~$2$, then $C$ has a simple singular point
at~$O$.) Then, pulling back, one arrives at
$s_0^3f_0^6\Gf=(s_0f_0^2p)^3+(s_0^2f_0^3q)^2$, where
$f_0\in\Gamma(\HS;\CO_\HS(F_0))$ is a section defining~$F_0$, and
the representation as in the statement is obtained by cancelling
$s_0^3f_0^6$.

Conversely, since~$S_0$ is contracted, any representation as in
the statement is pushed forward to a torus structure of~$C$.
\endproof

Consider a germ~$\Gf$ at a singular point~$P$ of type~$\bA_{3k-1}$
(respectively,~$\bE_6$) and fix local coordinates $(x,y)$ in which
$\Gf$ is given by $x^{3k}+y^2$ (respectively, $x^4+y^3$). Let, in
the same coordinate system, $p$ be a semiquasihomogeneous germ of
type $(k,1)$ (respectively, $(2,1)$), and let $q$ be adjacent to
a semiquasihomogeneous germ of type $([\frac12(3k+1)],1)$
(respectively, $(2,2)$).

\lemma\label{intersection}
In the notation above, assume that $\Gf-q^2=ph_1$ or
$\Gf-q^2=p^2h_2$.
Then $(h_1\cdot p)_P\ge\frac12(3k+1)$
\rom(\,$(h_1\cdot p)_P\ge3$ for $P$ of type~$\bE_6$\rom) or,
respectively, $(h_2\cdot p)_P\ge k$
\rom(\,$(h_2\cdot p)_P\ge2$ for $P$ of type~$\bE_6$\rom), where
$(a\cdot b)_P$ stands for the local intersection index at~$P$ of
the curves $\{a=0\}$ and $\{b=0\}$.
\endlemma

\proof
Under the assumptions, the Newton polygon of $\Gf-q^2$ at~$P$ is
contained in that of~$\Gf$, and a simple analysis gives a `lower
bound' for the Newton polygons of~$h_1$ and~$h_2$. If the singular
point~$P$ of~$\Gf$ is of type~$\bE_6$, then $h_1$ and~$h_2$ must
be adjacent to semiquasihomogeneous germs of type $(3,2)$ and
$(2,1)$, respectively. If $P$ is of type $\bA_{3k-1}$, then $h_2$
is adjacent to a semiquasihomogeneous germ of type $(k,1)$, and
the Newton polygon of~$h_1$ is a subset of the minimal Newton
polygon containing points $(0,2)$, $([\frac12(k+1)],1)$, and
$(2k,0)$. Now, the estimates for the local intersection indices
are immediate.
\endproof

\lemma\label{intersection.A2}
In the notation of Lemma~\ref{intersection}, assume that $P$ is of
type~$\bA_2$ and that $q$ is
only adjacent to a semiquasihomogeneous germ of type $(1,1)$. Then
$(h_1\cdot p)_P\ge1$ and $(h_2\cdot p)_P\ge1$.
\endlemma

\proof
The proof is similar to that of Lemma~\ref{intersection}. Now, we
can only assert that $\Gf-q^2$ is singular at~$P$ and, hence,
$h_1$ vanishes at~$P$. This proves the first estimate; for the
second one, we need to show that $h_2$ vanishes at~$P$ as well.
Assume the contrary. Then the singularity of~$\Gf$ at~$P$ is
equivalent to that of $p^2+q^2$, which is $\bA_{2r-1}$,
$r=(p\cdot q)_P$. This is a contradiction.
\endproof

\proof[\subsection{Proof of Theorem~\ref{main.J10}}]\label{s.proof}
Fix an irreducible $\bJ$-sextic~$C$ whose fundamental group
factors to~$\GD_6$, and consider its trigonal model~$B$.
The set of
singularities~$\Sigma_B$ is given by Corollary~\ref{D6}. If
$\Sigma_B=3\bA_2\splus\ldots$, we pick and fix three cusps and
ignore all other singular points of~$B$.

One has $\dim\ls|S_0+2F|=3$. Hence, there exists a curve
$Q=\{p=0\}\in\ls|S_0+2F|$ such that the germ of~$p$ at each
singular point~$P$ of~$B$ is as in Lemma~\ref{intersection}. This
curve is necessarily irreducible and, hence, nonsingular. (Indeed,
any reducible curve in $\ls|S_0+2F|$ is a union of~$S_0$ and two
fibers, and such a curve cannot pass through all singular points
of~$B$ with multiplicities prescribed above.)
In particular, from the B\'ezout theorem it follows that
the local intersection indices of~$Q$ and~$B$ at each point are
exactly as in the lemma, \ie, $(Q\cdot B)_P=2k$ if $P$ is of
type~$\bA_{3k-1}$, and $(Q\cdot B)_P=2$ if $P$ is of type~$\bE_6$.

Next, one has $\dim\ls|S_0+3F|=5$, and unless
$\Sigma_B=3\bA_2\splus\ldots$, there is a curve
$K=\{q=0\}\in\ls|S_0+3F|$ such that the germ of~$q$ at each
singular point~$P$ of~$B$ is as in Lemma~\ref{intersection}. If
$\Sigma_B=3\bA_2\splus\ldots$, we choose~$q$ as in
Lemma~\ref{intersection} at two of the three cusps and as in
Lemma~\ref{intersection.A2} at the third one. One can see that $K$
does not contain~$Q$ as a component; indeed, otherwise $K$ would
split into~$Q$ and a fiber~$F$, and a curve of this form cannot
satisfy all the conditions imposed. Since $Q\cdot K=3$, the
B\'ezout theorem implies that, at each singular point~$P$ of~$B$,
one has $(Q\cdot K)_P=k$ if $P$ is of
type~$\bA_{3k-1}$, and $(Q\cdot K)_P=2$ if $P$ is of type~$\bE_6$.

Let $\Gf\in\Gamma(\HS;\CO_\HS(3S_0+6F))$ be a section whose zero
set is~$B$.
Comparing the local intersection indices, one observes that the
restrictions $\Gf|_Q$ and $s_0q^2|_Q$ have the
same zero divisor. (We multiply~$q^2$ by~$s_0$ to make it a
section of the same line bundle as~$\Gf$; the restriction $s_0|_Q$
is a constant.) Hence, after an
appropriate rescaling, the restriction $(\Gf-s_0q^2)|_Q$ is
identically zero, \ie, $\Gf-s_0q^2=ph_1$ for some
$h_1\in\Gamma(\HS;\CO_\HS(2S_0+4F))$. Then
Lemmas~\ref{intersection} and~\ref{intersection.A2} imply that
$H_1\cdot Q\ge5$, where $H_1=\{h_1=0\}$. Hence, $H_1$ contains~$Q$
as a component, \ie, $\Gf-s_0q^2=p^2h_2$, and, applying
Lemmas~\ref{intersection} and~\ref{intersection.A2} and the
B\'ezout theorem once more, one concludes that the curve
$H_2=\{h_2=0\}$ coincides with~$Q$. Thus, after another rescaling,
$\Gf=p^3+s_0q^2$, and the statement of Theorem~\ref{main.J10} follows from
Lemma~\ref{reduction}.
\endproof

\Remark\label{rem.torus}
Theorem~\ref{main.J10} holds for reducible $\bJ$-sextics as well.
With the exception of the two degenerate families mentioned after
Lemma~\ref{simple}, the fundamental group of a reducible
$\bJ$-sextic~$C$ factors to~$\GD_6$ if and only if
the set of singularities of the trigonal model~$B$ of~$C$ is
$\Sigma_B=\bA_5\splus\bA_2\splus\bA_1$: it is the only root system
in~$\bE_8$ with $\bE_8/\Sigma_B$ having both $2$- and $3$-torsion.
(Formally, one would need to replace Corollary~\ref{dihedral} with an
analog of Theorem~3.5.1 in~\cite{degt.Oka}, which
would assert that, as in the irreducible case, the fundamental
group factors to~$\GD_6$ if and only if $\bE_8/\Sigma_B$ has
$3$-torsion.)
Then,
the proof given in Section~\ref{s.proof} extends literally to
$\Sigma_B=\bA_5\splus\bA_2\splus\bA_1$: one just
ignores the $\bA_1$-point.
\endRemark

\proof[\subsection{Proof of Theorem~\ref{main.J10}}]\label{s.proof'}
Unless the set of singularities~$\Sigma_B$ is $4\bA_2$,
there is a unique quotient
$\pi_2(\HS\sminus(B\cup S_0))\onto\GD_6$, see
Corollary~\ref{dihedral} and Lemma~\ref{Sigma.E8}, and
Section~\ref{s.proof} produces a unique torus structure. If
$\Sigma_B=4\bA_2$, in the construction of Section~\ref{s.proof}
one can choose any three of the four cusps; then the fourth one
is necessarily an outer singularity. Thus, one obtains four
distinct torus structures.
On the other hand, in
$\Tors(\CK_B,\ZZ/3\ZZ)=(\ZZ/3\ZZ)^3$ there are exactly four subgroups
isomorphic to~$\ZZ/3\ZZ$. The correspondence between the two sets
is established as in~\cite{degt.Oka}: a cusp~$P$ is inner (outer)
if and only if the composition
$\pi_1(U_P\sminus C)\to\pi_1(\Cp2\sminus C)\to\GD_6$ is
(respectively, is not) an
epimorphism (where $U_P$ is a Milnor ball about~$P$). In terms of
the cubic surface $V\subset\Cp3$ ramified at~$C$, the three inner
cusps are the cusps (Whitney pleats) of the projection $V\to\Cp2$,
whereas the outer one is the projection of an $\bA_2$-singular
point of~$V$.
\endproof

\section{Proof of Theorems~\ref{special}--\ref{abelian}}\label{S.groups}

We start with a description of two additional models of $\bJ$-sextics.

\Remark[The associated cubic]\label{s.cubic}
Let~$B$ be a trigonal model, and assume that $B$ has a triple
singular point~$P$. Blow it up, and blow down the fiber
through~$P$ and the exceptional section~$S_0$. The result
is~$\Cp2$, and the proper transform of~$B$ is a cubic
$\barD_3\subset\Cp2$.
The inverse transformation is determined by a point
$\barP\in\Cp2\sminus\barD_3$ (the image of~$S_0$) and a line $\barL$
through~$\barP$ (the image of the exceptional divisor over~$P$).
The triple $(\barD_3,\barP\in\barL)$ is called the
\emph{associated cubic} of~$B$.

The fiber of~$\HS$ through~$P$ contracts to a point; the other
fibers are in a one to one correspondence with the lines
through~$\barP$ other than~$\barL$.
\endRemark

\Remark[The associated quartic]\label{s.quartic}
This construction is similar to the previous one, but now we
start with a
double singular point~$P$ of~$B$. The proper transform of~$B$ is a
quartic curve $\barD_4\subset\Cp2$, and the inverse transformation
is determined by a nonsingular point $\barP\in\barD_4$. (The
exceptional divisor over~$P$ projects to the tangent~$\barL$ to~$\barD_4$
at~$\barP$.) The pair $(\barP\in\barD_4)$ is called the
\emph{associated quartic} of~$B$. The fibers of~$\HS$ other than
that through~$P$ transform to the lines through~$\barP$ other
than~$\barL$. Furthermore,
the birational map establishes a diffeomorphism
$\HS\sminus(B\cup S_0\cup F_P)\cong\Cp2\sminus(\barD_4\cup\barL)$,
where $F_P$ is the fiber of~$\HS$ through~$P$. As a consequence,
there is an epimorphism
$$
\pi_1(\Cp2\sminus(\barD_4\cup\barL))
 \onto\pi_1(\HS\sminus(B\cup S_0)).
$$
\endRemark

\lemma\label{special.fiber}
If a trigonal model~$B$ has a singular
fiber of type~$\tilde\bA_0^{**}$ or~$\tilde\bA_1^*$, then the
group $\pi_1(\HS\sminus(B\cup S_0))$ is abelian. If $B$ has a
singular fiber of type~$\tilde\bA_2^*$, then there is an epimorphism
$\GB_3\onto\pi_1(\HS\sminus(B\cup S_0))$.
\endlemma

\proof
The fundamental group $\pi_1(\HS\sminus(B\cup S_0))$ can be found
using van Kampen's method~\cite{vanKampen} applied to the ruling
of~$\HS$.
Remove a nonsingular fiber~$F_0$, pick another nonsingular
fiber~$F'$, and pick a generic section~$S$ disjoint from~$S_0$ and
from the critical points of the projection $B\to\Cp1$.
Let $G=\pi_1(F'\sminus(B\cup S_0),F'\cap S)$, and let $\Ga_1$,
$\Ga_2$, $\Ga_3$ be a standard set of generators of~$G$.
(Clearly, $F'\sminus(B\cup S_0)$ is a real plane with three
punctures.) Let $F_1,\ldots,F_r$ be the singular fibers of~$B$.
For each~$F_j$, dragging~$F'$ about~$F_j$ and keeping the base
point in~$S$ results in a certain automorphism $m_j\:G\to G$,
called the \emph{braid monodromy} about~$F_j$. Then, the group
$\pi_1(\HS\sminus(B\cup S_0),F'\cap S)$ has a representation of
the form
$$
\bigl<\Ga_1,\Ga_2,\Ga_3\bigm|
 \text{$m_j(\Ga_i)=\Ga_i$, $i=1,2,3$, $j=1,\ldots,r$, and
 $[\gamma_0]=1$}\bigr>,
$$
where $\gamma_0$ is a small circle in~$S$ about~$F_0$. (The class
$[\gamma_0]$ can be expressed in terms of~$\Ga_i$, but this
expression is irrelevant for our purposes.)

The same approach can be used to find the group
$\pi_1(U_F\sminus(B\cup S_0))$, where~$F$ is a fiber of~$\HS$,
singular or not, and $U_F$ is a tubular neighborhood of~$F$. The
resulting representation is
$\<\Ga_1,\Ga_2,\Ga_3,\,|\,\text{$m(\Ga_i)=\Ga_i$, $i=1,2,3$}>$,
where $m$ is the
local braid monodromy about~$F$.
An immediate and well known consequence
is the fact that the inclusion homomorphism
$\pi_1(U_F\sminus(B\cup S_0))\to\pi_1(\HS\sminus(B\cup S_0))$ is
onto.

The local monodromy can easily be found using model
equations;
for the fibers as in the statement, it is
given by the following expressions:
\widestnumber\item{$\tilde\bA_0^{**}:$}
\setbox0\hbox{$\tilde\bA_0^{**}:$}
\dimen0\wd0
\def\bx#1{\hbox to\dimen0{$#1$\hss}}
\roster
\item"\bx{\tilde\bA_0^{**}:}"
$\Ga_1\mapsto\Ga_2$,\quad $\Ga_2\mapsto\Ga_3$,\quad
$\Ga_3\mapsto\Pi\Ga_1\Pi^{-1}$,
\item"\bx{\tilde\bA_1^*:}"
$\Ga_1\mapsto\Ga_3$,\quad $\Ga_2\mapsto\Ga_3\Ga_2\Ga_3^{-1}$,\quad
$\Ga_3\mapsto\Pi\Ga_1\Pi^{-1}$,
\item"\bx{\tilde\bA_2^*:}"
$\Ga_1\mapsto\Ga_3$,\quad $\Ga_2\mapsto\Pi\Ga_1\Pi^{-1}$,\quad
$\Ga_3\mapsto\Pi\Ga_2\Pi^{-1}$,
\endroster
where $\Pi=\Ga_1\Ga_2\Ga_3$. Now, it is obvious that, in the first
two cases, the group $\pi_1(U_F\sminus(B\cup S_0))$ is free
abelian (of rank one and two, respectively), and in the last case,
$\pi_1(U_F\sminus(B\cup S_0))\cong\GB_3$.
\endproof

\corollary\label{E12}
If $C$ is a plane sextic with a singular point of type~$\bE_{12}$
or~$\bE_{13}$, then the fundamental group $\pi_1(\Cp2\sminus C)$
is abelian.
\endcorollary

\proof
Any sextic with a singular point of type~$\bE_{12}$ or~$\bE_{13}$
is a $\bJ$-sextic, and the statement follows from
Proposition~\ref{B=C}, Lemma~\ref{special.fiber}, and
Table~\ref{table1}.
\endproof

\lemma\label{triple.point}
If an irreducible trigonal model~$B$ has a triple point
and the set of singularities of~$B$ is not $\bE_6\splus\bA_2$,
then the group $\pi_1(\HS\sminus(B\cup S_0))$ is abelian.
\endlemma

\proof
Consider the associated cubic $(\barD_3,\barP\in\barL)$
and denote by~$\barP'$ the intersection of~$\barL$ and an
inflection tangent to~$\barD_3$ other than~$\barL$.
Moving~$\barP$ along~$\barL$ towards~$\barP'$ deforms~$B$ to a
trigonal model~$B'$ with a singular fiber of
type~$\tilde\bA_0^{**}$, and the statement of the lemma
follows from
Lemma~\ref{special.fiber}. The only exception is the case when
$\barL$ is the only inflection tangent
to~$\barD_3$; then $\barD_3$ has a cusp and
$\Sigma_B=\bE_6\splus\bA_2$.
\endproof

\lemma\label{simple.node}
If an irreducible trigonal model~$B$ has a simple node~$\bA_1$
and the set of singularities of~$B$ is not $3\bA_2\splus\bA_1$,
then the group $\pi_1(\HS\sminus(B\cup S_0))$ is abelian.
\endlemma

\proof
Let $(\barP\in\barD_4)$ be the associated quartic constructed
using a node~$P$. Unless $\barD_4$ is a three cuspidal quartic, it
has an inflection point~$\barP'$, see~\cite{quintics}, and
moving~$\barP$ towards~$\barP'$ deforms~$B$ to a trigonal model
with a singular fiber of type~$\tilde\bA_1^*$. Now, the statement
follows from Lemma~\ref{special.fiber}.
\endproof

\proof[\subsection{Proof of Theorem~\ref{special}}]\label{proof.special}
The sextics mentioned in the statement are precisely those whose
trigonal model has the set of singularities $2\bA_4$. Let~$C$ be
such a sextic. Denote by~$B$ its trigonal model and let
$G=\pi_1(\Cp2\sminus C)=\pi_1(\HS\sminus(B\cup S_0)$. Due to
Corollary~\ref{D6}, $G$ factors to~$\GD_{10}$; since also
$G/[G,G]=\ZZ/6\ZZ$, in fact $G$ must factor to
$\GD_{10}\times(\ZZ/3\ZZ)$.

Take for~$P$ one of the two singular points and
consider the associated quartic $(\barP\in\barD_4)$. The set of
singularities of~$\barD_4$ is $\bA_4\splus\bA_2$, and the
tangent~$\barL$ at~$\barP$ passes through the $\bA_2$-point. Such
pairs $(\barD_4,\barL)$ do exist, see~\cite{quintics}; this proves
the existence of sextics. Furthermore, there is an epimorphism
$\tilde G=\pi_1(\Cp2\sminus(\barD_4\cup\barL))\onto G$,
see Section~\ref{s.quartic}, and, according
to~\cite{groups}, \smash{$\tilde G$} is the semidirect
product given by the exact sequence
$$
1@>>>\FF_5[t]/(t+1)@>>>\tilde G@>>>\ZZ@>>>1,
$$
where $t$ is the conjugation by the generator of~$\ZZ$. The
largest quotient of~$\tilde G$ whose abelinization is $\ZZ/6\ZZ$
is again $\GD_{10}\times(\ZZ/3\ZZ)$. This completes the proof.
\endproof

\proof[\subsection{Proof of Theorem~\ref{abelian}}]\label{proof.abelian}
Let~$C$ be a sextic as in the statement, let~$B$ be its trigonal
model, and let $\Sigma_B\subset\bE_8$ be its set of singularities.
Then the group $\bE_8/\Sigma$ has neither $2$-torsion
(since $C$ is irreducible, Proposition~\ref{irreducible}),
nor $3$-torsion
(since $C$ is not of torus type, Theorem~\ref{main.J10}),
nor $5$-torsion (since $C$ is not as in Theorem~\ref{special});
thus, due to Lemma~\ref{E8}, the embedding
$\Sigma_B\subset\bE_8$ is primitive.

In view of Lemmas~\ref{triple.point} and~\ref{simple.node}, one
can assume that $B$ has neither triple points nor nodes. (The
exceptional cases in Lemmas~\ref{triple.point}
and~\ref{simple.node} are both of torus type.) Thus,
$\Sigma_B=\bigoplus\bA_{p_i}$, $p_i\ge2$, and applying Nori's
theorem~\cite{Nori} to the irreducible sextic with the set of
singularities $\bJ_{10}\splus\Sigma$ (\ie, considering a
non-singular distinguished fiber~$F_0$),
one rules out all
possibilities with $\sum(p_i+1)<9$. After this, the only set of
singularities left is $\Sigma_B=\bA_4\splus\bA_3$.

Let $\Sigma_B=\bA_4\splus\bA_3$. Take for~$P$ the $\bA_3$-point
and consider the associated quartic $(\barP\in\barD_4)$. The set
of singularities of~$\barD_4$ is $\bA_4\splus\bA_1$, and the
tangent~$\barL$ at~$\barP$ passes through the $\bA_1$-point.
According to~\cite{groups}, the group
$\pi_1(\Cp2\sminus(\barD_4\cup\barL))$ is abelian; hence, so is
$\pi_1(\Cp2\sminus C)$, see the epimorphism in
Section~\ref{s.quartic}.
\endproof

\proof[\subsection{Proof of Theorem~\ref{torus}}]\label{proof.torus}
According to Theorem~\ref{main.J10} and Corollary~\ref{D6}, an
irreducible
$\bJ$-sextic~$C$ satisfies the hypotheses of the theorem if and
only if its trigonal model $B$ has one of the following
sets of singularities:
$3\bA_2$, $3\bA_2\splus\bA_1$, $\bA_5\splus\bA_2$,
$\bA_8$, or $\bE_6\splus\bA_2$.

First, let us show that there is an epimorphism
$\GB_3\onto\pi_1(\HS\sminus(B\cup S_0))$. Due to
Lemma~\ref{special.fiber}, it is the case whenever $B$ has a
singular fiber of type~$\bA_2^*$. Otherwise, take for~$P$ one of
the $\bA_2$-points (or the only $\bA_8$-point) and consider the
associated quartic $(\barP\in\barD_4)$. Its set of singularities
is $\Sigma_B$ with one copy of~$\bA_2$ removed
(respectively, $\bA_6$), and the
tangent~$\barL$ at~$\barP$ is a double tangent (respectively,
passes through the $\bA_6$-point). In each case,
one has
$\pi_1(\Cp2\sminus(\barD_4\cup\barL))\cong\GB_3$,
see~\cite{groups},
and the desired
epimorphism is that of Section~\ref{s.quartic}.

One has $\GB_3/[\GB_3,\GB_3]=\ZZ$, and the central element
$\Delta^2\in\GB_3$ projects to $6\in\ZZ$. Thus, the largest
quotient of~$\GB_3$ whose abelinization is $\ZZ/6\ZZ$ is
$\GB_3/\Delta^2$, and one obtains epimorphisms
$$
\GB_3/\Delta^2\onto\pi_1(\HS\sminus(B\cup S_0))\onto\GB_3/\Delta^2.
$$
(The latter epimorphism is due to the fact that $C$ is of torus
type, see Theorem~\ref{th.main}.) On the other hand, the group
$\GB_3/\Delta^2=
 \operatorname{\text{\sl PSL}}(2,\ZZ)$ is Hopfian (as it is
 obviously residually finite). Hence, the two epimorphisms above
 are isomorphisms.
\endproof

\Remark
Alternatively, instead of referring to~\cite{groups}, one can
argue that, in all cases except $\Sigma_B=\bA_8$, the trigonal
model can be deformed into a curve with a singular fiber of
type~$\bA_2^*$. This procedure corresponds to deforming the
associated quartic~$\barD_4$ so that the double tangent~$\barL$
becomes a line intersecting~$\barD_4$ at a single point with
multiplicity four. Such a deformation exists due
to~\cite{quintics}.
\endRemark

\section{An application to the classification}\label{S.application}

The deformation classification of $\bJ$-sextics appeared
in~\cite{quartics}, but the proof has never been
published.
Most results of~\cite{quartics} related to
$\bJ$-sextics can be obtained by passing to the trigonal model~$B$ and
then, to the associated cubic or associated quartic. (As a matter
of fact, in
most cases the results are stated in terms of the associated cubic
or quartic, with a further reference to~\cite{quintics}.)
Most difficult are the cases when the $\bJ$-sextic has a
singular point of type~$\bJ_{2,1}$ or~$\bE_{12}$, so that one
needs to keep track of a singular fiber not passing through a
singular point of~$B$, see Table~\ref{table1}. In these cases, one
should add the transform~$\barL'$ of the singular fiber in
question to the associated quartic~$\barD_4$, thus reducing the
problem to the deformation classification of reducible sextics
$\barD_4+\barL+\barL'$ with simple singularities, see~\cite{JAG}.

Here, we present a simple proof of the `non-existence' statements
of~\cite{quartics}.

\proposition\label{no.J10}
There are no $\bJ$-sextics with the following sets of
singularities\rom: $\bJ_{2,i}\splus\bA_3\splus2\bA_2$,
$\bJ_{2,i}\splus\bA_4\splus2\bA_2$,
$\bJ_{2,i}\splus\bA_6\splus\bA_2$ \rom($i=0,1$\rom),
and $\bJ_{2,1}\splus4\bA_2$.
\endproposition

\proof
In the first three cases, the trigonal model~$B$ of the curve
would have set of singularities $\bA_3\splus2\bA_2$,
$\bA_4\splus2\bA_2$, or $\bA_6\splus\bA_2$. None of these lattices
admits an embedding to~$\bE_8$, see Corollary~\ref{embedding}. In
the last case, $B$ has four cusps, each cusp counting as a triple
singular fiber. On the other hand, $B$ always has twelve singular
fibers (counted with multiplicity); hence, in the case of four
cusps, there are no fibers of type~$\tilde\bA_0^*$, see
Table~\ref{table1}.
\endproof

\proposition\label{no.E12}
There are no $\bJ$-sextics with the sets of singularities
$\bE_{12}\splus\Sigma$, where either $\Sigma$
is one of the lattices listed in Lemma~\ref{Sigma.E8} or
$\Sigma=\bA_3\splus2\bA_2$, $\bA_4\splus2\bA_2$, or
$\bA_6\splus\bA_2$.
\endproposition

\proof
If $\Sigma=\Sigma_B$ is one of the lattices listed in
Lemma~\ref{Sigma.E8}, the fundamental group of the curve has a
dihedral quotient, see Corollary~\ref{D6}. This contradicts to
Corollary~\ref{E12}. The other three sets of singularities do not
admit an embedding to~$\bE_8$, see Corollary~\ref{embedding}.
\endproof

\Remark
For the sets of singularities $\Sigma=3\bA_2$, $\bA_5\oplus\bA_2$,
$\bA_8$, and $2\bA_4$, Proposition~\ref{no.E12} can be interpreted
as follows: a trigonal model~$B$ with one of these sets of
singularities cannot be deformed so that two simplest singular
fibers (of type~$\tilde\bA_0^*$) come together to form a singular
fiber of type~$\tilde\bA_0^{**}$. (For the other sets of
singularities listed in Lemma~\ref{Sigma.E8} this statement is
obvious as the curve has at most one type~$\tilde\bA_0^*$ singular
fiber.) Similarly, using the part of Corollary~\ref{E12}
concerning~$\bE_{13}$, one can see that, in the case
$\Sigma_B=3\bA_2\oplus\bA_1$, the $\bA_1$-point cannot join the
remaining singular fiber of type~$\tilde\bA_0^*$ to form a
singular fiber of type~$\tilde\bA_1^*$.
\endRemark

\enddocument